\documentclass[11pt,a4paper]{amsart}

\usepackage{latexsym}
\usepackage[dvips]{graphics}
\usepackage{color}
\usepackage{epsfig}

\usepackage{amssymb}
\usepackage{amsthm}
\usepackage{amsmath}
\usepackage{amstext}
\usepackage{amscd}
\usepackage{epic}

\theoremstyle{plain}
\newtheorem{thm}{Theorem}[section]
\newtheorem{prop}[thm]{Proposition}

\newtheorem{lem}[thm]{Lemma}

\newtheorem{theorem*}{Theorem}[]

\theoremstyle{definition}
\newtheorem{defn}[thm]{Definition}

\theoremstyle{remark}
\newtheorem{rem}[thm]{Remark}
\newtheorem{example}[thm]{Example}
\newtheorem{examples}[thm]{Examples}

\newcommand{\N}{\mathbb{N}}
\newcommand{\R}{\mathbb{R}}

\newcommand{\I}{\mathcal{I}}
\newcommand{\J}{\mathcal{J}}
\newcommand{\Oo}{\mathcal{O}}

\newcommand{\inv}{^{-1}} \def\({(\!(} 
\def\){)\!)}\DeclareMathOperator{\ord}{ord\, }
\DeclareMathOperator{\mult }{mult\,}

\DeclareMathOperator{\dist}{dist\,}

\date{13 March, 2009.}

\title{On the non-analyticity locus of 
an arc-analytic function}
\author{Krzysztof Kurdyka and Adam Parusi\'nski}
\address{
Laboratoire de Math\'ematiques, UMR 5175 du CNRS, Universit\'e de Savoie,
Campus Scientifique, 73 376 Le Bourget--du--Lac Cedex, France,
}
\email{kurdyka@univ-savoie.fr}

\address {Laboratoire Angevin de Recherche en Math\'ematiques, UMR 6093
du CNRS,  Universit\'e d'Angers,
   2, bd Lavoisier, 49045 Angers cedex, France}

\email{adam.parusinski@univ-angers.fr}

\thanks {This research was done at the 
Mathematisches Forschungsinstitut Oberwolfach during a
stay within the Research in Pairs Programme from
January 20 to February 2, 2008. We would like 
to thank the MFO for excellent working conditions.}

\begin{document}

\maketitle

 \begin{abstract}
 A function is called arc-analytic if it is real analytic on each real analytic arc.  In real analytic geometry there 
 are many examples of arc-analytic functions that are not real analytic.  Arc analytic functions appear while studying 
 the arc-symmetric sets and the blow-analytic equivalence.  
 In this paper we show that the non-analyticity locus of an arc-analytic function is arc-symmetric. 
 We discuss also the behavior of the non-analyticity locus under blowings-up.    By 
 a result of  Bierstone and Milman a big class of arc-analytic function, namely those that satisfy a polynomial equation with 
 real analytic coefficients, can be made analytic by a sequence of global blowings-up with smooth centers.  
 We show that these centers can be chosen, at each stage of the resolution,  inside  the non-analyticity 
 locus.  
 \end{abstract}



\bigskip
\section {Introduction. }
\label{Intro}

Let $X$ be a real analytic manifold.
A function $f:X\to \R$ is called
\emph{arc-analytic},   cf. \cite{kurdyka1}, if for every real analytic
$\gamma:(-1,1)\to X$ the composition $f\circ \gamma$ is
analytic.
The arc-analytic functions are closely related to blow-analytic functions of Kuo, cf. 
\cite{kuo2}.  In particular,
we have the following result, conjectured for the functions with semi-algebraic graphs in 
\cite{kurdyka1},
and shown in \cite{bm1}.

\begin{thm}\label{bm}
Let $X$ be a nonsingular real analytic manifold and let $f: X\to \R$ be an 
arc-analytic function on $X$.  Suppose that 
\begin{equation*}\label{algebraic}
G(x,f(x)) =0,
\end{equation*}
where
$$
G(x,y) = \sum_{i=0}^p g_i(x) y^{p-i}
$$
is a nonzero polynomial in $y$ with coefficients $g_i(x)$ which are analytic  functions 
on $X$.
Then there is a mapping $\pi:X'\to X$ which is a composite of a locally finite sequence 
of blowings-up with nonsingular closed centers, such that $f\circ \pi$ is  analytic.
\end{thm}

 Let $f:X\to \R$ be an arc-analytic subanalytic function.
 In this paper we study \emph{the set $S(f)$ of non-analyticity
  of $f$}.  By definition, $S(f)$ is the complement of the set $R(f)$
  of points $p\in X$, such that $f$ as a germ is real analytic at $p$.   It is known (cf.
\cite{tamm}, \cite{kksmooth}, \cite{BM88}) that $S(f)$ is
  closed and subanalytic.  It follows from \cite{bm1}  or \cite{parusinski}, that $\dim 
S(f) \le \dim X -2$.  As we show in
  Theorem \ref{arc-symmetric} below, $S(f)$ is arc-symmetric in the sense of 
\cite{kurdyka1}.    Theorem \ref{arc-symmetric} is shown in section \ref{sectionas}.

  We also study how the set of non-analyticity behaves under blowings-up with smooth 
centers.  This
  depends on whether the center is entirely contained  in $S(f)$ or not.  If it is not 
then
  the non-analyticity lifts to the entire fiber, see Proposition \ref{blowing-up}.  Note 
that
 Theorem \ref{bm} can be also derived from \cite{parusinski}.  Using the method of 
\cite{parusinski} and
 Proposition \ref{blowing-up}  we show
  the following refinement of
 Theorem \ref{bm}.

 \begin{thm}\label{extra}
In Theorem \ref{bm} we may require that the mapping $\pi:X'\to X$, that is
a locally finite composite $\pi = \cdots \circ \pi_k \circ \cdots \circ \pi_0$ of 
blowings-up with smooth centers,
satisfies additionally:

\smallskip
for every $k$ the center of $\pi_{k+1}$ is contained in the locus of non-analyticity of 
$f\circ \pi_0 \circ \cdots
\circ \pi_k$.
 \end{thm}
 \bigskip


\subsection{Algebraic case}
Theorem \ref{bm} can be stated in the real algebraic version, see \cite {bm1}.  In this 
case if  we assume that
$X$ is a nonsingular real algebraic variety and that the coefficients
$g_i$ are regular then we may require that  $\pi$ is a finite composite of blowings-up 
with nonsingular
algebraic centers.

In the algebraic case we cannot require that the centers of blowings-up  are entirely 
contained in
the non-analyticity loci as Example \ref{contre} shows.

An analytic function on $X$  is called
\emph{Nash} if its graph is semialgebraic.  It is called
\emph{blow-Nash} if it can be made Nash after composing with a
finite sequence of blowing-ups with smooth nowhere dense
regular centers.  Thus the algebraic version of Theorem \ref{bm}, cf. \cite{bm1},  says 
that the function with semi-algebraic graph is arc-analytic if and only if it is
blow-Nash.  Nash morphisms and manifolds form a natural category that contains the 
algebraic one, cf. \cite{BCR}.
 We note that our refinement of the statement of Theorem \ref{bm} holds in the Nash 
category.

\begin{thm}\label{Nashbm}
Let $X$ be a  Nash manifold and let $f:X\to \R$ be an arc-analytic function
on $X$.  Suppose that
\begin{equation*}
G(x,f(x)) =0,
\end{equation*}
where
$$
G(x,y) = \sum_{i=0}^p g_i(x) y^{p-i}
$$
is a nonzero polynomial in $y$ with coefficients $g_i(x)$ which are Nash  functions on 
$X$.
Then there is a finite composite
 $\pi = \cdots \circ \pi_k \circ \cdots \circ \pi_0$ of blowings-up of nonsingular Nash 
submanifolds,  such that
 for every $k$ the center of $\pi_{k+1}$ is contained in the locus of non-analyticity of 
$f\circ \pi_0 \circ \cdots
\circ \pi_k$, and  $f\circ \pi$ is  Nash.
\end{thm}


\subsection{Subanalytic case}

Less is known for an arc-analytic function with subanalytic graph if it does not
 satisfy an equation \eqref{algebraic}.  It is known that an arc-analytic subanalytic
 function has to be continuous and can be made real analytic by composing with
finitely many local blowings-up with smooth centers, see  \cite{bm1}  or
\cite{parusinski} (we refer the reader to these papers for a precise
statement).  It is not known whether these blowings-up can be made global
that is whether the arc-analytic subanalytic functions coincide with the
family of blow-analytic functions  of T.-C.Kuo, see e.g.
\cite{kuo2}, \cite{fukuikoikekuo}, \cite{fukuipaunescu}.    It is also  not known,
whether the centers of such  blowings-up can be chosen in the locus of
non-analyticity of the function.

We present below in Example \ref{nasty} a subanalytic arc-analytic function that cannot 
be made analytic,  even locally,  by  a blowing-up of a coherent ideal.  In particular, it cannot satisfy 
an equation of type \eqref{algebraic}.

\medskip
\subsection{Examples}\label{examples}
\begin{examples}\label{manyexamples}
The function  $f:\R^2 \to \R$ , $f(x,y) =\frac{x^3}{x^2+y^2} $ for
$(x,y) \ne (0,0)$ and $f(0,0) =0$, is arc-analytic but not differentiable at the origin.

The function $g(x,y) = \sqrt {x^4 + y^4}$ is arc-analytic but not $C^2$.
This example is due to E. Bierstone and P.D. Milman.

The function  $h:\R^2 \to \R$ , $h(x,y) =    \frac{xy^5}{x^4+y^6}$ for  $(x,y) \ne (0,0)$ and $h(0,0) =0$ 
is arc-analytic but not lipschitz. This  example is due to L. Paunescu.

We generalize the first example as follows.
 Fix a real analytic Riemannian metric on $X$  and let $Y$ be a nonsingular real analytic 
subset of $X$.   Then $ d_Y^2: X\to \R$, the square of the distance to  $Y$, is a real analytic function 
on $X$.   Suppose that $Y$ is of codimension $\ge 2$ in $X$ and let $f:X\to \R$
 be an analytic function vanishing on $Y$ and not divisible by $d^2$. Then, $\frac {f^3} 
{d^2}$   vanishes on $Y$, is arc-analytic and   not analytic at the points of $Y$.  
Note that $\frac {f^3} {d^2}$ composed with the blowing-up of $Y$ is analytic.
 \end{examples}

 \begin{example}\label{contre}
  Let $ g(x,y) = y^2 +  x(x-1) (x-2)(x-3)$.  Then $g\inv (0) \subset \R^2$ is
irreducible and has two connected compact components,  denoted by $X_1$ and $X_2$.  These
connected components  that can be separated by
$h(x,y) = x-1.5$, that is $h<0$ on $X_1$ and $h>0$ on $X_2$.    For $\varepsilon>0$ 
sufficiently small,  $h^2 + \varepsilon g$ is strictly positive on $\R^2$.   Define
$$
g_1 (x,y) = \sqrt{h^2 +  \varepsilon g} + h.
$$
Then $g_1$ is  analytic, $0$ is a regular value of $g_1$ and $g_1\inv (0) = X_1$.   Moreover,
$g_1$ is Nash.  Then $f: \R^3 \to \R$ defined by
$$
f(x,y,z) = \frac{z^3} {z^2 + g_1^2(x,y)}
$$
for $(x,y,z)\ne 0$ and $f(0)=0$, is arc-analytic and $S(f) = X_1\times \{0\}$.  The function $f$ 
becomes analytic after blowing-up of  $S(f)$.
\end{example}

\begin{example} \label{nasty}
Let $\pi_0 : \widetilde \R^3 \to \R^3$ be the blowing-up of the origin and let $E$ be the 
exceptional divisor of $\pi_0$.
Let $C\subset E$ be a transcendental (the smallest algebraic subset of $E$ that contains 
$C$ is $E$ itself)   non-singular analytic curve and let $\pi_C: M\to \widetilde \R^3$ be
the blowing-up of $C$.   Let $f$ be an arc-analytic function on $\R^3$ such that the set 
of non-analyticity   of  $f\circ \pi_0$ is $ C$ and
$f\circ \pi_0\circ \pi_C$ is analytic.  Such a function can be constructed as follows.  
Using the last remark of
Examples \ref{manyexamples} we may construct an arc-analytic function $g:\widetilde \R^3\to 
\R$ such that $S(g) = C$.  Then we may set  $f(x,y,z)= (x^2+y^2+z^2) \, g( \pi_0\inv (x,y,z))$.

Such $f$, as a germ at $0$,
cannot be made analytic by a single blowing-up of an ideal.
Indeed,  suppose contrary to our claim that there exists an ideal $\I$ of 
$\R\{x_1,x_2,x_3\}$ such that $f\circ \pi_{\I}$ is analytic,
where $\pi_{\I}$ denotes the blowing-up of $\I$.  Multiplying $\I$ by the maximal ideal 
at $0$ we may assume that $\pi_{\I}$  factors through $\pi_0$, i.e. $\pi_{\I}= \pi_{\J}\circ \pi_0$, 
where $\J$ is a sheaf of coherent ideals centered
on an algebraic subset $Y$ of $E$.  We may assume that $\dim Y\le 1$. Thus the blowing-up 
of $\J$,  $\pi_{\J} :  M_{\J} \to   \widetilde R^3 $ is an isomorphism over the complement of $Y$ 
that contradicts the  construction of $f$.
\end{example}


\newpage
\section {Arc-meromorphic mappings. }
\label{arc-mero}

In this section  \emph{subanalytic}  mean  subanalytic at infinity. Let us 
recall, \cite{tamm}, \cite{kksmooth}, that a subset  $A$ of
$\R^n$ is called \emph {subanalytic at infinity} if  $A$ is subanalytic in some algebraic 
compactification of $\R^n$. (Then in  fact it is subanalytic in every algebraic compactification of $\R^n$.)
All functions and mappings are supposed to be subanalytic, that is their graphs are 
subanalytic at infinity.

\begin{defn}\label{arcmerodef}
Let $U$ be an open subanalytic  subset of $\R^n$.  An everywhere defined
subanalytic mapping $f:U\to \R^m$ is 
called
\emph{arc-meromorphic} if  for any analytic arc
$\gamma :(-1,1) \to U$ there exists  a discrete set $D\subset (-1,1)$ and
$\varphi$ an meromorphic function on $(-1,1)$ with poles contained in $D$ and such that
 $f\circ \gamma = \varphi$  on $(-1,1)\setminus D$.
  Note that it may happen that  $f\circ \gamma$ does not coincide with $\varphi$ at some 
points of $D$ and may be
  at these points discontinuous.
\end{defn}

\begin{example} 
The function  $f:\R^2 \to \R$ defined by  $f(x,y) =\frac{xy}{x^2+y^2}$ for
$(x,y) \ne (0,0)$ can be extended to an arc-meromorphic function on $\R^2$ by assigning 
any value at the 
origin.  Then it becomes discontinuous at $(0,0)$  even if for every analytic arc $\gamma :(-1,1) \to \R^2$, 
$\gamma(0)=(0,0)$, $f\circ \gamma$ extends to  an analytic function.
\end{example}

\begin{rem}\label{contmero}
If $f$ is  an arc-meromorphic and continuous function on an open connected
set $U\subset \R^n$, then $f$ is arc-analytic.
\end{rem}

\begin{rem}\label{identity}
Let $f$ and $g$ be arc-meromorphic functions on an open connected
set of $U$.  Assume that  $f=g$ on an open non-empty subset $U\subset \R^n$,
 then $f=g$  except on a nowhere dense subanalytic subset of $U$.
\end{rem}

\begin{lem}\label{lojasinqlem}
Let $U$ be an open bounded subanalytic subset in $\R^n$ and $f:U\to \R^m$
be an arc-meromorphic mapping. Then there exists $\Gamma \subset \R^n$ a closed nowhere dense 
subanalytic set, $N\in \N$ and $C>0$ such that
\begin{equation}\label{lojasinq1}
|f(x)|\le C \dist(x,\Gamma)^{-N}, \, x\in U\setminus \Gamma.
\end{equation}
In particular we can take as $U$ the complement of the non-analyticity locus of  $f$.
\end{lem}

\begin{proof}
It is well-known (cf. e.g. \cite{hironaka}, \cite{lojfourier})  that there exists a stratification of
$\R^n$ which is compatible with $\overline U$ and  such that $f$ is analytic on each 
stratum contained in $U$.
We take as $\Gamma $ the union of all strata  contained in $\overline U$ of dimension 
less than $n$. Let  us consider the function defined as follows:  $g(x) = |f(x)|$ if $|f(x)|  \le 1$,
and $g(x) = |f(x)|^{-1} $ if $|f(x)|\ge 1$.
Then  $h(x):= dist(x,\Gamma)g(x)$ is a subanalytic and continuous function on $\overline U$ 
which is compact. Moreover, if $\dist(x,\Gamma)=0$ then $h(x)=0$. Therefore, by the classical
\L ojasiewicz's inequality (cf. e.g. \cite{hironaka}, \cite{BM88}) for subanalytic functions, there exist $N\in 
\N$  and $c>0$ such that
\begin{equation}
h(x) \ge c \dist (x,\Gamma)^{N+1}, \, x\in U.
\end{equation}
Thus inequality \eqref{lojasinq1}   follows with $C= \max\{1/c,M\}$, where $  M=\sup_{x\in U} 
dist(x,\Gamma)^N $.
\end{proof}

We state now an  auxiliary lemma on arc-meromorphic functions  in two variables.

\begin{lem}\label{killsinglem}
Let $U$ be an open subanalytic subset in $\R^2$ and let  $f:U\to \R^m$
be an arc-meromorphic mapping.  Then for  any $a\in U$ there exists
a neighborhood  $V$ of $a$ and an analytic function $\varphi : V\to \R$, 
$\varphi\not\equiv 0$, such that $\varphi f $  is arc-analytic.
\end{lem}

\begin{proof} Let $\Gamma$   be the subanalytic set  associated  to
$f$  by  Lemma  \ref{lojasinqlem}. Clearly we may assume that $a\in \Gamma$,
otherwise  $f$ is analytic at $a$ and the statement is trivial. Since $\dim \Gamma = 
1$, by a result of {\L}ojasiewicz's  \cite{ihes} (see also \cite{klz}), the set $\Gamma$ is actually 
semianalytic. Then there exists a neighborhood  $V$ of $a$ and an analytic function $\psi : V'\to \R$, 
$\psi\not\equiv 0$,   which vanishes on $V'\cap\Gamma$. Hence  for some compact neighborhood
$V\subset V'$ of $a$ there exists $c>0$ such that  
$$
|\psi(x)| \le c \dist(x, \Gamma),\, x\in V.
$$
 (This is a consequence of the main value theorem). Put $\varphi=\psi^{N+1}$, then by  Lemma
 \ref{lojasinqlem} the function $\varphi f$ is continuous on $V$. Clearly $\varphi f$ is 
 arc-meromorphic,  so by Remark \ref{contmero} this function is arc-analytic.
\end{proof}

\begin{prop}\label{arcmeroderiv}
Let $f:U\to \R$ be an arc-meromorphic function, where $U$ is an open subset in $\R^n$.
Assume that $f$ is analytic with respect to the variable $x_1$. Then the 
function $\frac{\partial f}{\partial x_1}:U\to \R$ is again arc-meromorphic.
\end{prop}

\begin{proof}
First observe that by \cite{kksmooth}  the function $\frac{\partial f}{\partial x_1}$
is (globally) subanalytic. To prove that  $\frac{\partial f}{\partial x_1}$ is 
arc-meromorphic
let us fix an analytic arc $\gamma : (-1,1) \to U$. We define an arc-meromorphic function
$g: V\to \R$ by
$g(s,t)= f(\gamma(t) + se_1)$, where $e_1 = (1, 0, \dots 0)$ and $ V$ is an open
neighborhood of $\{0\} \times (-1,1)$ in $\R^2$. Clearly
$$
\frac{\partial f}{\partial x_1} (\gamma (t))= \frac{\partial g}{\partial s} (0,t).
$$
 By  Lemma \ref{killsinglem} there exist a neighborhood   $V$ of $(0,0)$ and  an analytic function   
$\varphi: V\to \R$ such that $h:= \varphi g $ is  arc-analytic on $V$.  Since  $\dim S(h)\le 0$,
for any $t\ne 0$ sufficiently small  $h$ is analytic at $(0,t)$, but of course also 
$\varphi$ is analytic at $(0,t)$.
Since $ g(s)$  is analytic with respect to $s$ it follows that  $g= h/\varphi$ is 
actually analytic at $(0,t)$  for  any
$t\ne 0$ sufficiently small.
By \cite{bm1} there exists a map  $\pi:M\to \R^2$, which is a finite composition of 
blowing-ups of points, such that $h\circ \pi $ is analytic.
Consider the arc $\eta (t): = (0,t)$ and let $\tilde \eta (t) \in M$ be the unique analytic arc such that
$\pi\circ \tilde \eta = \eta$.
The chain rule gives
\begin{equation}
d_{\tilde \eta(t)} h\circ \pi =  (d_{ \eta(t)} h) \circ (d_{\tilde \eta(t)}  \pi).
\end{equation}
Note that $d_{\tilde \eta(t)}  \pi$ is invertible for $t\ne 0$, moreover the map  $t\mapsto (d_{\tilde \eta(t)}  \pi)^{-1}$ 
is meromorphic.  It follows that $t\mapsto d_{ \eta(t)} h $ is meromorphic. In particular
$t\mapsto \frac{\partial h}{\partial s}(0,t)$ is meromorphic. We have
$$
\frac{\partial h}{\partial s}(0,t) = \varphi \frac{\partial g}{\partial s}(0,t) + 
g\frac{\partial \varphi}{\partial s}(0,t).
$$
Since $\varphi(0,t)\ne 0$ for $t\ne 0$, 
the map  $t\mapsto \frac{\partial g}{\partial s}(0,t)$ is meromorphic and Proposition 
\ref{arcmeroderiv}  follows.
\end{proof}

\begin{rem}\label{allderivatives}
A repeated application of Proposition \ref{arcmeroderiv} shows that for every $k\in \N$, 
 $$
 \frac{\partial^k f}{\partial x_1^k}:U\to \R
 $$
 is arc-meromorphic.  Moreover, there exists a subanalytic stratification 
 $\mathcal S$ of $U$ such that for every stratum $S\in \mathcal S$ and every $x\in S$ there is $\varepsilon >0$ and  
 a neighborhood $V$ of $x$ in $S$ such 
 that $f(x+se_1)$ is an analytic function of $(x,s) \in V\times (-\varepsilon,\varepsilon)$.  In particular, for every  $k\in \N, $ ${\partial^k f} /{\partial x_1^k}:U\to \R$  is analytic on the strata of $\mathcal S$.  
 \end{rem}
\section{The non-analyticity locus of an arc-analytic function is arc-symmetric.}\label{sectionas} 
Let $U\subset\R^n$ be open and let $f:U\to \R$  be arc-analytic with subanalytic graph.
We  denote by $S(f)$ the non-analyticity set of $f$ and by $R(f)$ its complement in $U$.  
Then  $S(f)$  is closed in $U$ and by \cite{tamm} (see also  \cite{kksmooth},
\cite{bm1}) it is a subanalytic set.  It follows from \cite{bm1}  or \cite{parusinski} that $\dim S(f) \le n -2$.

\begin{thm}\label{arc-symmetric}
Let $\gamma : (-\varepsilon , \varepsilon ) \to U$  be an analytic arc such that $\gamma 
(t) \in R(f)$ for $t<0$.  Then
$\gamma (t) \in R(f)$ for $t>0$ and small.  In other words, $S(f)$ is arc-symmetric 
subanalytic in the sense of  \cite{kurdyka1}.
\end{thm}

For the proof we need some basic properties of Gateaux differentials.
For each $k\in \N$ we consider
\begin{equation}\label{hfunctions}
h_k(x,v) = \frac{1}{k!} \partial^k_v f (x) = \frac 1 {k!} \frac {d^k}{dt^k} 
f(x+tv)_{|t=0}. 
\end{equation}

\begin{prop}\label{dirdercor}
Let $f:U\to \R$ be an arc-analytic function.  Then for  any  $k\in \N$ the function
$h_k (x,v) : U\times \R^n \to\R$  is arc-meromorphic.
\end{prop}

\begin{proof}
Let $(x(t), v(t))$ be an analytic arc in $ U\times \R^n$.  Define an arc-analytic function
$g(s,t)= f(x(t) + s v(t))$.  Then
$$
h_k (x(t) ,v(t) ={\frac 1{k!} } \frac {\partial^k } {\partial s^k}  g(t,s)_{|s=0}
$$
that is meromorphic  by Proposition \ref{arcmeroderiv}.
\end{proof}

For  $x\in U$, $k\in \N$ we denote 
$$
h_{x,k} (v) = h_k (x,v) = \frac{1}{k!} \partial^k_v f (x)
$$
Note that  $h_{x,k}$ is $k$-homogeneous function. If $f$ is analytic at $x$,  then $h_{x,k}$ is polynomial.
We have  also the inverse which is Bochnak-Siciak Theorem, see \cite {BS2},
 which states that if $h_{x,k}$ is polynomial for each $k\in \N$, then $f$ is
 analytic at $x$. Traditionally if  $h_{x,k}$ is polynomial then it is called
the Gateaux  differential of $f$ at $x$ of order $k$.
 
We call $h_{x,k}$ {\it generically polynomial} if it is equal to a polynomial except
on a nowhere dense subanalytic (and homogenous) subset of $\R^n$.
Note that, by Remark \ref{identity},  $h_{x,k}$  is generically polynomial if it coincides 
with a polynomial  on an open nonempty set.

\begin{prop}\label{Gateauxgeneric}
Let $f:U\to \R$ be an arc-analytic function, where $U$ is an open subset in $\R^n$.
Let $\gamma : (-\varepsilon,\varepsilon) \to U$ be an analytic arc and $k\in \N$.
If $h_{\gamma (t) ,k}$ is generically polynomial for 
$t \in (-\varepsilon,0)$, 
then  there exists a finite set $F_k \subset (0, \varepsilon)$
 such that $h_{\gamma (t) ,k}$ is generically   polynomial
 for each $t \in (0, \varepsilon)\setminus F_k$.
\end{prop}

\begin{proof}
Let  $\R_k[x_1,\dots,x_n]$ denote the space of  homogenous polynomials
of degree $k$ and let $d_k=\binom {n+k-1}{n}$ denote  its dimension.
We need the classical multivariate interpolation.  

\begin{lem}\label{mvinterpol} 
 There exists an algebraic nowhere dense subset $\Delta \subset (\R^n)^{d_k} $ such that 
for $V=(v^1, \ldots , v^{d_k})\in  (\R^n)^{d_k} \setminus \Delta$ the map $\Psi_V :\R_k[x_1,\dots,x_n] \to  \R^{d(k)} $ 
given by 
$$
\Psi_{V} (P) = (P(v^1), \ldots , P(v^{d_k})). 
$$
is a linear isomorphism.  
 \qed \end{lem}

Fix $V=(v^1, \ldots , v^{d_k})\in  (\R^n)^{d_k} \setminus \Delta$ generic and denote  
$\Phi_V = \Psi_V\inv  :\R^{d(k)} \to \R_k[x_1,\dots,x_n]$.  We define an arc-meromorphic map 
$P_k :(-\varepsilon,\varepsilon) \to 
\R_k[x_1,\dots,x_n]$
by
$$
P_k (t) : =  \Phi_V (h_k(\gamma(t), v^1), \dots, h_k(\gamma(t), v^{d(k)}) ).  
$$
The map $p_k;(-\varepsilon,\varepsilon)\times\R^n \to \R$, where $p_k(t,v) =P_k(t)(v)$ is arc-meromorphic.  
If $V$ is sufficiently generic then, for $t\in (-\varepsilon,0)\setminus \{\text{finite set}\}$, 
$p_k(t)$ coincides with 
$h_{\gamma(t),k}$.  Since they both are arc-meromorphic, by Remark \ref{identity}
they coincide  on $(-\varepsilon,\varepsilon)\times\R^n \setminus Z_k$, where $Z_k$ is a 
closed subanalytic set with $\dim Z_k \le n$. Hence there  exists a finite set 
$F_k \subset (0, \varepsilon)$  such that  for $t \in (0, \varepsilon)\setminus F_k$ the intersection 
$Z_k\cap (\{t\}\times \R^n)$ is of  dimension less than $n$. Thus,  for each 
$t \in (0, \varepsilon)\setminus F_k$ the  function  $h_{\gamma (t) ,k}$ is  {generically} polynomial,
as claimed.
\end{proof}


The following  proposition is a version of the mentioned above Bochnak-Siciak Theorem, \cite{BS2}.

\begin{prop}\label{bsgeneric} 
If for every $k$ there is a nonempty open subset $V_k\subset \R^n$ and a homogeneous polynomial 
$P_k$ of degree $k$ such that $h_{x,k}\equiv P_k$ on $V_k$,  then $f$ is analytic at $x$.  
\end{prop} 

\begin{proof}
We first show that $\sum_k P_k(v)$ is convergent in a neighborhood of $0\in \R^n$.

We may assume that $x$ is the origin.  Let $\pi_0$ be the blowing up of the origin,
$\pi_0(y,s) = (sy,s), s\in \R,\, y\in\R^{n-1},$ in a chart.   The function $\tilde f 
(y,s):= f(\pi(y,s))$, defined in  a neighborhood $U'$
of the exceptional divisor $E:s=0$, is arc-analytic.   The set of non-analyticity  of 
$\tilde f$, denoted by $\tilde S$, is  closed subanalytic and
of codimension at least $2$.   For $y\notin \tilde S$, $\tilde f$ is analytic in a 
neighborhood of $(0,y)$ and,  moreover, by analytic continuation,
\begin{equation}\label{razem}
h_{x,k}(v) = P_k(v) \qquad \text { for} \quad v=t(y,1),  \,  t\in \R, \, y \notin \tilde S.
\end{equation}

Fix $A'$ an open non-empty  subset of  $E$ such that the closure of $A'$ does not 
intersect $\tilde S$.
Let $A\subset \R^n$ be the cone over $A'$.  Then, by \eqref{razem}, $\sum_k P_k(v)$ is 
convergent in any compact subset of $A$.  The convergence in  a neighborhood of $0$ in $\R^n$
follows from the following lemma.

\begin{lem}\label{alibaba}
Let $V\subset \R^n$ be starlike with respect to the origin, $a\in V$, and suppose that 
$$
|P_k(v)|\le L \quad \text { on } V'=a+V.
$$
Then
$$
|P_k(v)|\le L \quad \text { on } \frac 1 {2e} V.
$$
\end{lem}

\begin{proof}
Since $P_k$ is homogeneous of degree $k$
\begin{equation}\label{binoms}
P_k(v)= \frac 1 {k!} \sum_{s=0}^{s=k} (-1)^{k-s} \binom k s P_k(a+sv).
\end{equation}
Indeed, \eqref{binoms} can be shown recursively on $k$ using Euler's formula as follows.
First note \eqref{binoms} holds for $a=0$ and the derivative of the RHS of \eqref{binoms} 
with respect to $a$  equals
\begin{equation}\label{binoms2}
0 = \frac 1 {k!} \sum_{s=0}^{s=k} (-1)^{k-s} \binom k s Q(a+sv),
\end{equation}
where $Q(x) = \sum _{i-1}^n a_i \frac{\partial P_k}{\partial x_i} (x)$ is a homogeneous 
polynomial of degree $k-1$.  By the inductive assumption
\begin{eqnarray*}
& &  \sum_{s=0}^{s=k} (-1)^{k-s} \binom k s Q(a+sv)   = \sum_{s=0}^{s=k-1} (-1)^{k-1-s} 
\binom {k-1} s Q(a+sv) + \\
&  & \quad  + \sum_{s=1}^{s=k} (-1)^{k-s} \binom {k-1}
 {s-1} Q(a+sv)  = -Q(v) + Q(v) = 0
\end{eqnarray*}
This shows \eqref{binoms}.
Thus, if $v\in \frac 1 k V$, $|P_k(v)|\le \frac 1 {k!} L \sum _{s=0} ^{k} \binom k s =  L  \frac  {2^k}{k!}$,
that means that for $v\in \frac 1 {2e} V$
$$
|P_k(v)|\le L \frac {(2k)^k} {k!} \frac 1 {(2e)^k}   \le L  .
$$
This ends the proof of lemma \ref{alibaba}.
\end{proof}Then  $\sum_k P_k(v)$ is an analytic function in a neighborhood of the origin that coincides with 
 $f$ on a set with non-empty interior.  Hence
$f(v) =\sum_k P_k(v)$ in a neighborhood of the origin. This shows proposition 
\ref{bsgeneric}.
\end{proof}

\begin{proof}[Proof of theorem \ref{arc-symmetric}]
We may assume that $\gamma$ is injective otherwise the image of $t>0$ equals the image of 
$t<0$ and the statement is obvious. Let $F:= \bigcup F_k$, where $F_k$ are finite subsets of
$(0, \varepsilon)$  given by Proposition  \ref{Gateauxgeneric}.  
Clearly the complement of $F$ is dense in $(0, \varepsilon)$, so 
by Proposition \ref{bsgeneric} our function $f$ is analytic  at $\gamma(t)$ for $t\in G$, 
where $G$ is an open dense subset of $(0, \varepsilon)$.  
Hence  theorem \ref{arc-symmetric} follows.
\end{proof}

Consider the subanalytic sets
\begin{align*}
\tilde R_{k_0} (f) =  \{x\in U; \, & \forall k\le k_0, \, h_{x,k}   \text { is   
generically polynomial }\}  ,  \\
  R_{k_0} (f) =  \{x\in U; \, & \forall k\le k_0, \, h_{x,k}   \text { is   
  polynomial }\}  .
\end{align*}

 Clearly  $\tilde R_{k+1}(f)\subset \tilde R_k(f)$ and 
$ R_{k+1}(f)\subset  R_k(f)$.  We recall from \cite{kksmooth}  the following  result 

\begin{prop}\label{kkgateaux}{\rm [ \cite{kksmooth}, Proposition 4.4] }
Let $f:U\to \R$ be a subanalytic (not necessarily arc-analytic) function on an open bounded $U\subset  \R^n$.
Then for any compact $K\subset U$  there is $k\in \N$ such that $R(f) \cap K =  R_k (f) \cap K$.
\end{prop}

\begin{prop}\label{bsgenericfinite}
For any compact $K\subset U$ there is $k\in \N$ such that $R(f) \cap K =  \tilde R_k (f) \cap K$.
\end{prop}

\begin{proof}
By Remark \ref{allderivatives} there exists a stratification $\mathcal S$ of $U\times S^{n-1}$ such that for every $k$, 
$h_k$ is analytic on the strata.  Refining the stratification, if necessary, we may suppose that for every stratum $S\subset U\times  S^{n-1}$ its projection to $U$ has all fibers of the same dimension.  In the proof we use only these 
strata for which all the fibers of projection to $U$ are of maximal dimension $n-1$.  We denote the collection of them 
by $\mathcal S_n$ and their union as $Z$.    Now it is easy to adapt the proof of  Lemma 6.1  of 
\cite{kksmooth} (based on multivariate interpolation) and show the following lemma. 

\begin{lem}
There are analytic subanalytic functions 
$$
w_i : U \times S^{n-1} \to \R, \quad i\in \N,
$$
analytic on each stratum of $\mathcal S$ such that $h_{x,i}$ is generically polynomial if and only if $w_i \equiv 0$ generically on 
$\{x\}\times S^{n-1}$.  \qed
	\end{lem}

Now Proposition \ref{bsgenericfinite} follows from Lemma 2.5 of 
\cite{kksmooth}   that shows that for every stratum there exist $k$ such that 
$$
\bigcap_{i=1}^\infty \{w_i=0\}= \bigcap_{i=1}^k \{w_i=0\}.
$$
\end{proof}

We complete this section with two results, one that controls the change of non-analyticity locus by blowings-up.   
This result will be crucial in the next section.  The last result of this section, Proposition \ref{1parameter}, though not used in this paper,  indicates a possible analogy between our approach and the theory of complex analytic functions. 

\begin{prop}\label{blowing-up}
Let $T=\{x_k=x_{k+1}=\cdots=x_n=0\}$ and let $\pi_T$ be the blowing-up of $T$.
Suppose that the origin is in the closure of $R(f)\cap T$ and that $f\circ \pi_T$ is 
analytic at least at one point of
 $\pi_T \inv (0)$ (hence on a neighborhood of this point).
 Then $f$ is analytic at $0$.
\end{prop}

\begin{proof}
 Let $\Pi: \R^n \times \R \times \R^n\to \R^n$ be given by $\Pi(x,t,v) =  x+tv$ and let
 $\Pi_T: T\times \R \times \R^{n}\to \R^n $ be the restriction of $\Pi$.   First we show 
that if  $f\circ \Pi_T$ is analytic at some points of  $\Pi_T\inv (0)\cap \{t=0\}$ and $0$ is in 
the closure of $R(f) \cap T$ then $f$ is analytic at $0$.   Indeed,
 suppose that $A'\subset \R^n$ has non-empty interior and suppose that $f\circ \Pi_T$ is 
analytic in a neighborhood  $\{0\}\times \{0\} \times A'$.  Let $h_k(x,v)$, $x\in T, v\in \R^n$, be defined by 
\eqref{hfunctions}.  Then $h_k$ is arc-meromorphic  and analytic on $A=U'\times A'$, where 
$U'$ is a small neighborhood of $0$ in $T$.   For each $k$, we define by Lemma \ref{mvinterpol}, 
 \begin{equation}\label{interpol}
 P_k (x,v)= \Psi_V\inv (h_k (x, v^1), \dots, h_k(x, v^{d(k)}) )(v) , 
 \end{equation}
where $v^1,\ldots, v^{d_k}\in  A'$ are generic.   Each $P_k$ is
 analytic on $A$ and equals $h_k$ for $x\in R(f)\cap T$.  Therefore $h_k (0,v)=P_k(0,v)$ 
for $v\in A'$ and the   claim follows from proposition \ref{bsgeneric}. 

 Thus it remains to show that $f\circ \Pi_T$ is analytic at some points of $\Pi_T\inv 
(0)\cap \{t=0\}$.   For this we factor $\Pi_T$ restricted to   $\{v_n\ne 0\}$ through $\pi_T$ and use the 
assumption on $\pi_T$.  Write $\pi_T$ in an affine chart $\pi_T (\tilde x,  y,s)
= (\tilde x,   s y, s)$, where  $\tilde x =(\tilde x_1, \ldots, \tilde x_{k-1})$,  $ y =   ( y_{k}, \ldots,y_{n-1})$ and $s\in \R$.
Then on these charts $\Pi_T= \pi_T\circ \varphi$, where
$$
(\tilde x,y,s) = \varphi (x,t,v) = (x+tv', \frac 1 {v_n} v'', tv_n) ,
$$
 where
$v'=(v_1,\ldots, v_{k-1})$, $v''=(v_k,\ldots, v_{n-1})$.
 Restricted to $t=0$, $\varphi$ is a surjective projection $(x,v)\to (x,\frac 1 
{v_n}v'')$ onto $s=0$.  Hence
 $R(f\circ \Pi_T) \cap \Pi_T\inv (0)\cap \{t=0\}\supset \varphi\inv (R(f\circ \pi_T)\cap 
\pi\inv (0))$ is non-empty.
 \end{proof}

\begin{prop}\label{1parameter}
Let $x=(x_1,x')\in \R\times \R^{n-1}$ and suppose that for every $x_1>0$ and small, 
$f(x_1,x')$ is analytic at $(x_1,0)$
as a function of $x'$.  Moreover, suppose that for $x_1>0$ and small we have a uniform 
bound
$$
|h_k ((x_1,0), v')| \le c^k ,   \qquad \text { for } \|v'\|\le \varepsilon, k\in \N,
$$
where $v'=(v_2, \ldots,v_n)$.
Then $f$ is analytic at the origin.
\end{prop}

\begin{proof}
The function $h_k ((x_1,0), v')$ is arc-meromorphic as a function of $x_1,v'$.
Moreover, since continuous arc-meromorphic
functions of one variable are analytic, using polynomial interpolation lemma, Lemma \ref{mvinterpol},  
 we may show that each  $h_k ((x_1,0),v')$ extends
to an analytic function $\Psi(x_1,v')$ defined in a neighborhood of $(0,0)$, such that for each $x$, 
$v' \to \Psi(x_1,v')$ is a homogeneous polynomial in $v'$.  Moreover, for $x_1>0$ and $\|x'\|< \varepsilon/c$
$$
f(x_1,x')= \sum _k h_k ((x_1,0), x')
$$
and the series on the right-hand side is convergent.

Fix any $k\in \N$ and $\|v'\|< \varepsilon/c$.  Then for $v=(1,v')$, $0<t<1$,
$$
f(tv)= \sum_{j=0}^\infty h_j ((t,0), tv') = \sum_{j=0}^\infty t^j h_j ((t,0), v') =  
\sum_{j=0}^k t^j  h_j ((t,0), v') + \varphi (t,v') ,
$$
where $\varphi$ is subanalytic and $O(t^{k+1})$.  Therefore for such $v$
\begin{equation}
H_k(0,v) := \frac 1 {k!} \frac {d^k}{dt^k} f(tv)_{|t=0} =  \frac 1 {k!} \frac {d^k}{dt^k} 
 \sum_{j=0}^k  h_j ((t,0), tv')_{|t=0}.
\end{equation}
Note that the right-hand side, and hence $H_k (0,v)$ as well, is a polynomial in $v$.  
Indeed, this follows from the fact that $x\to \sum_{j=0}^k  h_j ((x_1,0), x')$ 
is an  analytic function of $x$ and $H_k(0,v) $ coincides with its Gateaux differential.  Thus 
proposition \ref{1parameter} follows from proposition \ref{bsgeneric}.
\end{proof}



\medskip
\section {Proof of Theorem \ref{extra}. }
\label{locus}

 We may suppose that $U$ is connected.  We suppose also that  the coefficients $g_0$ and 
$g_p$ of $G$
 and {the discriminant $\Delta (x)$ of $G$ are not identically equal to zero.  By 
the resolution of singularities  \cite{hironakares}, \cite{bm2},
 \cite{wlodarczyk},   there is a locally finite sequence of blowings-up $\pi:U'\to U$ 
with nonsingular centers such
 that $(g_0 g_p \Delta)\circ \pi $ is normal crossings.  Thus Theorem \ref{bm} follows from 
the following.

 \begin{prop}\label{old}
 Let an  arc-analytic function $f(x)$ satisfy the equation \eqref{algebraic} with analytic 
coefficients $g_i$.
 If $g_0$, $g_p$ and $\Delta (x)$ are simultaneously normal crossings (and hence not 
identically equal to zero) then
 $f$ is real analytic.
 \end{prop}

Proposition \ref{old} was proven in \cite{parusinski} under an additional assumption 
$g_0\equiv1$, see the proof of Theorem 3.1 of \cite{parusinski}.  It is easy to reduce the proof to this case 
by replacing $f$ by $g_p f$.  Then, an argument of \cite{parusinski} shows that locally $f$ can be 
expand as a fractional power series.  Finally, an arc-analytic fractional power series is analytic, see the proof 
of Theorem 3.1 of \cite{parusinski}.  If the discriminant  of $G$ vanishes identically then we replace it by the first 
non-vanishing higher order discriminant.

 To show Theorem \ref{extra} we follow, for the product $h(x) = g_0(x) g_p(x)  \Delta (x) 
$,   the
 monomialisation procedure    of W{\l}odarczyk or Bierstone-Milman. In this procedure the 
centre of blowing-up
    is defined as a the locus of points where a local invariant is maximal.   Thus  
suppose that we have the following data described in a local system of coordinates $x_1, 
\cdots, x_n$ at the
 origin.   The function $h \circ  \pi$, where $\pi =  \pi_k \circ \cdots \circ \pi_0$, is 
of the form $h\circ \pi=x^A h_k$, where $h_k$ is the controlled transform by the 
preceding blowings-up.
 Let $m = \ord_x  h_k$.   We
 may assume that $H=\{x_n=0\}$ is a hypersurface of maximal contact.
 Then, using the notation $x=(x',x_n)$,
 $$
 h_k(x) = x_k^m + \sum _{j=0}^{m-2} c_j (x')  x_k^j ,
 $$
 and $\mult_0 c_i\ge m-i$.

     Let $C$ be  the next centre given by the procedure and denote by $\pi_C$ the
     blowing-up of $C$.       We show that it cannot happen that
 $0\in S(h\circ \pi )$ and $0\in \overline{ C\setminus  S( f\circ \pi)}$.
 Suppose, contrary to our claim, that this is the case.  Then, by Proposition 
\ref{blowing-up}, the fibre over the origin
  of  the blowing-up $\pi_C=\pi_{k+1}$ of $C$ is contained in $S( f\circ \pi \circ 
\pi_C)$.   Since $C$ is contained in the equimultiplicity locus of $h_k$,
  at the generic point $\pi_C\inv (0)$ the strict transform of $h_k$ is nonzero, and 
hence $h\circ \pi\circ \pi_C$
  is normal crossing.  This contradicts Proposition \ref{old}.

  Let $C'$ denote the connected component of $C$ containing $0$. Then either $C'\subset 
S( h\circ \pi)$ or
  $C'\cap S( h\circ \pi)= \emptyset$.   Thus Theorem \ref{extra} proven.
 \qed




\end{document}